\documentclass[12pt,a4paper]{article}
\usepackage{amsmath, amstext, amsxtra, amsthm, amscd, amsgen, amsbsy, amsopn,
amsfonts, latexsym, amssymb, amscd, array}

\def\c{\text{{\bf c}}}
\def\spc{Spin^{\c}}

\def\proofend{\hbox to 1em{\hss}\hfill $\blacksquare $\bigskip }
\def\powser#1{\lbrack \lbrack #1 \rbrack \rbrack }

\newtheorem{theorem}{Theorem}[section]
\newtheorem{proposition}[theorem]{Proposition}

\newtheorem{remarks}[theorem]{Remarks}

\newtheorem{corollary}[theorem]{Corollary}

\newtheorem{conjecture}[theorem]{Conjecture}

\def\Z{{\mathbb Z}}
\def\R{{\mathbb R}}
\def\Q{{\mathbb Q}}
\def\C{{\mathbb C}}

\def\b{bun\-dle}
\def\pb{principal \b }
\def\vb{vector \b }

\def\mfd{manifold}

\def\twoindex#1#2#3{\underset {#1}{\overset {#2}{#3}}}

\def\oha{\mathcal H}

\def\paperref#1#2#3#4#5#6{\text{#1:} #2, {\em #3} {\bf#4} (#5)#6}
\def\bookref#1#2#3#4#5#6{\text{#1:} {\em #2}, #3 #4 #5#6}
\def\preprintref#1#2#3#4#5{\text{#1:} #2, {\em #3}, #4 (#5)}

\def\LFF{Lefschetz fixed point formula}
\begin{document}

\title{Homotopy complex projective spaces with $Pin(2)$-action}
\author{Anand Dessai}
\date{}
\maketitle

\begin{abstract}
\noindent
Let $M$ be a manifold homotopy equivalent to the complex projective 
space $\C P^m$. Petrie conjectured that $M$ has standard total Pontrjagin 
class if $M$ admits a non-trivial action by $S^1$. We prove the 
conjecture for $m<12$ under the assumption that the action extends to 
a nice $Pin(2)$-action with fixed point. The proof involves 
equivariant index theory for $Spin^c$-manifolds and Jacobi functions 
as well as classical results from the theory of transformation groups.
\footnote{keywords: homotopy complex 
projective spaces, rigidity, transformation groups, equivariant index theory, Jacobi functions}
\end{abstract}

\section{Introduction}
Let $M$ be a smooth closed manifold homotopy equivalent to the complex 
projective space $\C P^m$ and let $G$ be a compact Lie group which 
acts smoothly and non-trivially on $\C P^m$. We consider the problem 
to determine how close $M$ and $\C P^m$ are as differentiable 
manifolds if $M$ also supports a non-trivial action by $G$. In this 
paper we give an answer to this problem if $G$ is equal to $Pin(2)$, 
the normalizer of a maximal torus in $S^3$, and the dimension of $M$ 
is less than $24$. 

By simply-connected surgery theory one knows that for fixed $m\geq 3$ 
the set of diffeomorphism classes of homotopy $\C P^m$'s is infinite 
and partitioned into finite subsets by their total Pontrjagin class. 
In view of this classification one may think of a homotopy complex 
projective space $M$ as being close to $\C P^m$ if the total 
Pontrjagin class of $M$ is standard, i.e. if $p(M)$ takes the standard 
form $(1+x^2)^{m+1}$, where $x$ is a generator of $H^2(M;\Z )$. 

To be more specific on the problem above we ask the following strong 
(resp. weak) question: Is the total Pontrjagin class of a homotopy complex projective space $M$ standard 
(resp. standard up to finite ambiguity), if $M$ also supports a 
non-trivial action by $G$? 

The strong question has been answered in special cases by Petrie, 
Hattori, Masuda and many others (cf. \cite{Do} for a survey). Some of 
their work is stated below. Here we point out that by a result of 
Petrie (cf. \cite{Pe2}) the total Pontrjagin class of $M$ is standard 
if $M$ admits an effective action by the $m$-dimensional torus. The 
weak question has been considered for certain $S^3$-actions in 
\cite{DBF}. The questions above are motivated by the following 
conjecture of Petrie (cf. \cite{Pe1}).

\bigskip 
\begin{conjecture} Let $M$ be a homotopy $\C P^m$. If $M$ supports
 a non-trivial smooth $S^1$-action then its total Pontrjagin class is standard. 
\end{conjecture}

Petrie's conjecture has been verified in several special cases (again 
cf. \cite{Do} for a survey). In particular, it holds for $m\leq 4$ 
(cf. \cite{De}, \cite{Ja}) or if the number of fixed point components 
of the $S^1$-action is $\leq 4$ (cf. \cite{Wa}, \cite{Yo}, 
\cite{TsWa}, \cite{Ma}). In contrast the conclusion of the conjecture 
is known to fail if the $S^1$-action is only locally linear or if the 
group acting on $M$ is finite of arbitrary size (cf. \cite{DoMa}, 
\cite{Do}). 

One may argue that the positive results towards Petrie's conjecture 
described above merely express the general principle that the presence 
of a symmetry group imposes strong restrictions on the 
topology if the dimension of the symmetry group is `large' 
compared to the dimension of the manifold or if the orbit structure is 
`simple'. Note however, that by a result of Hattori (cf. \cite{Ha}) 
the conclusion of Petrie's conjecture holds in {\it any} dimension if 
$M$ admits an $S^1$-equivariant stable almost complex structure with 
standard first Chern class. Hattori's argument is based on vanishing 
results for equivariant $Spin 
^c$-Dirac operators derived from the 
\LFF \ in equivariant index theory (cf. \cite{AtSe}, \cite{AtSiIII}). 

In \cite{DTOP} generalizations of the rigidity theorems for elliptic 
genera (cf. \cite{Wi}, \cite{BoTa}, \cite{Hi}, \cite{Li}) to 
$Spin^c$-manifolds were given. These results led to a proof of 
Petrie's conjecture (cf. \cite{DAN}) in the special case that the 
first Pontrjagin class of $M$ is standard and $M$ carries a smooth $Pin(2)$-action
 with fixed point which is cohomologically trivial and 
almost effective (i.e. the action has finite kernel). We shall call 
such an action {\bf nice with fixed point}. Note that such actions 
exist on any $\C P^m$ for $m>1$. The main purpose of this paper is to 
show that in small dimensions the assumption on the first Pontrjagin 
class can be removed. This gives a partial answer to the strong 
question stated above. 

\bigskip
\begin{theorem}\label{main} Let $M$ be a smooth closed homotopy $\C
P^{m}$ which supports a nice $Pin(2)$-action with fixed point. If 
$m<12$ then the total Pontrjagin class of $M$ is standard. 
\end{theorem}

To prove the theorem we first give a lower bound for the first 
Pontrjagin class of $M$ (see Section \ref{linear bound}) which we 
obtain by relating equivariant index theory to the theory of Jacobi 
functions (see Section \ref{jacobi}). Then we combine results of 
\cite{DAN} with the homotopy invariance of $p_1(M)\bmod 24$ to 
complete the proof (see Section \ref{proof}). 

Note that for $m$ odd the $Pin(2)$-action on $M$ is trivial on 
cohomology if and only if $Pin(2)$ acts by orientation preserving 
diffeomorphisms. For $m$ even the $Pin(2)$-action has a fixed point if 
it acts trivially on cohomology (see Section \ref{equiroots}). In 
particular, a non-trivial smooth $S^3$-action on a homotopy $\C 
P^{2N}$ induces a nice $Pin(2)$-action with fixed point by restricting 
the $S^3$-action to the normalizer of a maximal torus. Hence, Theorem 
\ref{main} implies

\bigskip 
\begin{corollary}\label{maincor} Let $M$ be a smooth closed homotopy $\C
P^{2N}$ which supports a non-trivial smooth $S^3$-action. If $2N<12$ 
then the total Pontrjagin class of $M$ is standard.\proofend
\end{corollary}

This paper is structured in the following way. In the next section we 
give some information on the weights of equivariant vector bundles 
with vanishing first Pontrjagin class. In Section \ref{jacobi} we 
define certain series of $S^1$-equivariant twisted $Spin^c$-Dirac 
operators and express their indices in terms of Jacobi functions. In 
Section \ref{linear bound} we establish a lower bound for the first 
Pontrjagin class of a cohomology $\C P^m$ which admits a nice 
$Pin(2)$-action with fixed point. In the final section we give the 
proof for a slightly more general version of Theorem \ref{main}. 

\section{Weights and the first Pontrjagin class}\label{equiroots}
In this section we give some information on the weights of 
$Pin(2)$-equivariant vector bundles with vanishing first Pontrjagin 
class. Let $M$ be a $2m$-dimen\-sional smooth oriented closed manifold 
with smooth $Pin(2)$-action. We assume that the $Pin(2)$-action is 
trivial on integral cohomology. This property guaranties that the 
action lifts to complex line bundles over $M$. For the induced action 
of $S^1\subset Pin(2)$ let $Y$ denote a connected component of the 
fixed point manifold $M^{S^1}$. 

Let $\hat \xi \to M$ be a $Pin(2)$-equivariant $s$-dimensional complex 
\vb . Since $Y$ is a trivial $S^1$-space the restriction of $\hat 
\xi $ to $Y$ (viewed as an $S^1$-equivariant \vb ) splits as a finite 
sum $\hat 
\xi _{\vert Y}=\sum _{k\in \Z }\hat \xi _k\otimes \lambda ^k$. Here $\hat \xi _k$ is 
a complex 
\vb 
\ over $Y$ which is trivial as an $S^1$-space and $\lambda $ denotes 
the standard complex one-dimensional representation of $S^1$ (to 
lighten the notation we suppress the dependence of $\hat \xi _k$ on $Y$). Let $u_1, 
\ldots ,u_s$ denote the roots of $\sum _k\hat \xi _k $ defined using 
the splitting principle. Then the equivariant roots of $\hat \xi 
_{\vert Y}$ are defined as $u_1 +\omega _1\cdot z,\ldots ,u_s+\omega 
_s\cdot z$, where $\omega _i$ is equal to $k$ if $u_i$ is a root of 
$\hat \xi _k $ and $z$ is a formal variable. We call $\omega 
_1,\ldots ,\omega 
_s$ the weights of $\hat \xi $ at $Y$. Note that the 
character of the complex 
$S^1$-representation $\hat \xi 
_{\vert y}$, $y\in Y$, is equal to $\sum _i\lambda ^{\omega _i}$.

Next assume $\hat \xi $ is the complexification of an oriented 
$Pin(2)$-equivariant real $2t$-dimensional \vb \ $\xi $. Then $\hat 
\xi $ is invariant under conjugation and the equivariant roots of 
$\hat \xi 
_{\vert Y}$ occur in pairs $(u_i +\omega _i\cdot z, u_{t+i} +\omega _{t+i}\cdot z)$, 
$i=1,\ldots ,t$, where $u_{t+i} +\omega _{t+i}\cdot z=-(u_i +\omega 
_i\cdot z)$. We call $\pm (u_i+\omega _i\cdot z)$, $i=1,\ldots ,t$, 
the equivariant roots and $\pm \omega _i $ the weights of $\xi _{\vert 
Y}$. A spectral sequence argument for the Borel construction of $\xi $ 
shows (cf. \cite{DAN}, Prop. 3.7)

\bigskip
\begin{proposition}\label{negind} If the first Pontrjagin class $p_1(\xi )$ is torsion then 
$\sum _{i=1}^t \omega _i^2$ is independent of $Y\subset M^{S^1}$.\proofend 
\end{proposition} 

We apply the proposition in the case that $M$ is a cohomology $\C 
P^m$, i.e. $H^*(M;\Z )\cong H^*(\C P^m;\Z )$, and $p_1(M)$ is equal to $-n\cdot x^2$,
 where $x\in H^2(M;\Z )$ is a generator and $n$ is a non-negative 
integer. Let $\pm (x_i+m_{Y,i}\cdot z)$, $i=1,\ldots ,m$, denote the 
equivariant roots of $TM_{\vert Y}$. Consider the complex line bundle 
$\gamma $ over $M$ with $c_1(\gamma)=x$. Since $Pin(2)$ acts trivially on integral 
cohomology there exists a unique lift of the $Pin(2)$-action to $\gamma $ 
(this follows from \cite{HaYo}, cf. \cite{DAN}, Prop. 3.6). Let $a_Y$ be the weight of the induced 
$S^1$-action on $\gamma $ at $Y$. 

Note that by the assumption on $p_1(M)$ the first Pontrjagin class of 
the bundle $\xi 
=TM +n\cdot 
\gamma $ vanishes. Hence, Proposition \ref{negind} implies

\bigskip
\begin{corollary} There is a constant $C\in \Z$ such that $\sum _{i=1}^m m_{Y,i}^2+n\cdot a_Y^2=C$
for any connected component $Y$ of $M^{S^1}$.\proofend 
\end{corollary}

Next we assume that the $Pin(2)$-action is nice with fixed point, i.e. 
we assume in addition that the $Pin(2)$-action has a fixed point 
$pt\in M$. A simple application of the Lefschetz fixed point formula 
for the Euler characteristic shows that such a fixed point exists if 
the Euler characteristic of $M$ is odd (cf. \cite{DAN}, Lemma 3.8). 
Let $Y_0$ denote the connected component of $M^{S^1}$ which contains 
$pt$. Note that the representation $\gamma 
_{\vert pt}$ is trivial since it is a complex one-dimensional 
$Pin(2)$-representation. Hence $a_{Y_0}$ vanishes and we conclude from 
the corollary above that 
\begin{equation}\tag{$\ast $}\sum _{i=1}^m 
m_{Y,i}^2+n\cdot a_Y^2=\sum _{i=1}^m m_{Y_0,i}^2\end{equation} for any 
fixed point component $Y$. This formula is used in Section \ref{linear 
bound} to give a lower bound for the first Pontrjagin class. 

\section{Twisted $\mathbf {Spin^c}$-Dirac operators and Jacobi functions}\label{jacobi}
In this section we consider certain series of equivariant twisted 
$Spin^c$-Dirac operators closely related to elliptic genera 
and describe their indices in terms of Jacobi functions. Let $M$ be a 
$2m$-dimensional closed connected manifold with $Spin^c$-structure 
given by a $Spin^c(2m)$-\pb 
\ $P\to M$.\footnote{Together with a fixed isomorphism between the 
induced $SO(2m)$-\pb \ and the frame bundle.} The $Spin^c$-structure 
induces a complex line bundle over $M$ and we denote its first Chern class by $c\in 
H^2(M;\Z )$. 

Let $V\to M$ be a complex 
\vb \ of dimension $s$. We fix connections on $V$ and the $U(1)$-part of 
$P$. Let $\partial _{\c}\otimes V$ denote the associated twisted 
$Spin^c$-Dirac operator acting on sections of the tensor product of 
the complex spinor bundle and $V$. By the Atiyah-Singer index theorem 
(cf. \cite{AtSiIII}) its index $ind(\partial 
_{\c}\otimes V)$ is a topological invariant given by
$$ind(\partial _{\c}\otimes V)=\langle e^{c/2}\cdot \hat {\mathcal A}(M)\cdot ch(V),\mu _M\rangle
.$$
Here $\hat {\mathcal A}(M)$ denotes the multiplicative series for $M$ 
associated to the $\hat A$-genus, $\mu _M$ is the fundamental cycle of 
$M$ and $\langle \quad , \quad \rangle $ denotes the pairing between 
cohomology and homology. 

Next assume $M$ carries an $S^1$-action and the action lifts to the 
$Spin^c$-structure $P$ and to the complex vector bundle $V$. In this 
situation the twisted $Spin^c$-Dirac operator $\partial 
_{\c}\otimes V$ refines to an $S^1$-equivariant operator and its index 
refines to an element $ind_{S^1}(\partial _{\c}\otimes V)$ of the 
complex representation ring $R(S^1)$. For any topological generator 
$\lambda 
_0\in S^1$ the equivariant index 
$ind_{S^1}(\partial _{\c}\otimes V)(\lambda _0)$ may be computed from 
the 
\LFF 
\ (cf. \cite{AtSe}, \cite{AtSiIII}) in terms of local data at the 
fixed points
$$ind_{S^1}(\partial _{\c}\otimes V)(\lambda _0)=\sum _Y\tilde \nu _Y(\lambda _0).$$
Here the sum runs over the connected components $Y$ of the fixed point 
manifold $M^{S^1}$. To describe the local data $\tilde \nu _Y$ it is 
convenient to replace the $S^1$-action by the two-fold 
action.\footnote{This is not necessary but makes formulas easier. In particular, some 
functions which are only well defined on the covering $\C ^*\to 
\C ^*, 
\lambda \mapsto \lambda 
^2$, are defined on the base after passing to the two-fold action.} 
Having done so it follows from the \LFF \ that each local contribution 
$\tilde 
\nu 
_Y(\lambda )$ is a rational function in $\lambda \in \C $ which 
only depends on the restriction of the $Spin^c$-structure $P$ and the 
\b \ $V$ to $Y$. Since $ind_{S^1}(\partial _{\c}\otimes V)\in R(S^1)$ 
is a finite Laurent polynomial in $\lambda $ the sum $\sum _Y\tilde 
\nu _Y(\lambda )$ extends to a meromorphic function on $\C $ without 
poles on $\C 
^*$. 

Below we shall consider a certain series of $S^1$-equivariant twisted 
$Spin^c$-Dirac operators for which the equivariant index is related to 
a Jacobi function (see Prop. \ref{local data}). In the next section we 
employ this relation to study the first Pontrjagin class of a 
homotopy complex projective space which admits a nice $Pin(2)$-action 
with fixed point. 

We digress and recall the definition of Jacobi functions. Let $SL_2(\Z )$ act on $\Z^2$ by
 matrix multiplication from the 
right, i.e. $(\alpha ,\beta )\mapsto (\alpha ,\beta )A$ for $A\in SL_2(\Z )$, and let $\oha $ denote 
the upper half-plane.

A meromorphic function $F(\tau ,z)$ on $\oha 
\times 
\C $ is called a Jacobi function for $SL_2(\Z )\ltimes \Z ^2$ of weight $k$ and index $I$ if 
$$F(\tau ,z+\alpha \tau +\beta)=F(\tau ,z)\cdot e^{-2\pi i\cdot I\cdot (\alpha ^2\cdot \tau +2\alpha
z)}$$ for $(\alpha ,\beta )\in \Z ^2$ and
$$F(\frac {a\tau +b}{c\tau +d},\frac z{c\tau +d})=F(\tau ,z)\cdot (c\tau +d)^k\cdot e^{2\pi i\cdot I\cdot
\frac {c\cdot z^2}{c\tau +d}}$$ for $\left
(\begin{smallmatrix}a&b\\c&d\end{smallmatrix}\right )\in SL_2(\Z )$. In 
view of these equations one may also define Jacobi functions of weight $k$ and index $I$ as fixed 
points under an action of $SL_2(\Z )\ltimes \Z ^2$ on the ring of meromorphic functions
 on $\oha \times \C $ (cf.
\cite{EiZa} where the definition also involves conditions for the 
cusps). 

In topology Jacobi functions occur naturally as local contributions in 
the \LFF \ of elliptic genera for $Spin$, stable almost complex or 
$BO\langle 8\rangle $-manifolds (cf. \cite{Wi}, \cite{BoTa}, 
\cite{Hi}, \cite{Li}). 

We shall now consider some generalizations of elliptic genera to $S^1$-equi\-variant $\spc 
$-mani\-folds. As before let $V$ be an $S^1$-equivariant $s$-dimensional 
complex \vb \ over $M$. We define a $q$-power series ${\mathcal U}_V\in 
K_{S^1}(M)\powser q$ of virtual $S^1$-equivariant \vb s by 
$${\mathcal U}_V:=\twoindex {n=1} \infty \bigotimes S_{q^n}(\widetilde {TM}\otimes _\R \C )\otimes
\Lambda _{-1}(V^*)\otimes \twoindex {n=1} \infty \bigotimes \Lambda _{-q^n}(\widetilde
{V}\otimes _\R \C ).$$ Here $q$ is a formal variable,
$\widetilde E$ denotes the reduced \vb \ $E-\dim (E)$ and
$\Lambda _t:=\sum \Lambda ^i\cdot t^i$ (resp. $S_t:=\sum S^i\cdot t^i$) denotes the exterior
(resp. symmetric) power operation. The tensor product is, if not indicated otherwise, taken
over the complex numbers.

The index of the equivariant $Spin^c$-Dirac operator twisted with
${\mathcal U}_V$ is a $q$-power series of representations 
$ind_{S^1}(\partial _{\c}\otimes {\mathcal
U}_V)\in R(S^1)
\powser q$. By the \LFF \ the equivariant 
index at a topological generator $\lambda _0$ of $S^1$ is a sum of 
local data 
$$ind_{S^1}(\partial _{\c}\otimes {\mathcal U}_V)(\lambda _0)=\sum _Y\tilde \nu _Y(q,\lambda _0).$$
Each local datum $\tilde \nu _Y(q,\lambda )$ is an element of $\C
(\lambda )\lbrack \lbrack q\rbrack \rbrack $ which only depends on the 
restriction of the $Spin^c$-structure $P$ and ${\mathcal U}_V$ to $Y$. In 
order to explain its relation to Jacobi functions we need to introduce 
some notation for the local data at the $S^1$-fixed points. 

For a connected component $Y$ of $M^{S^1}$ define $d(Y):=\dim (Y)/2$ 
(since $M$ is of even dimension the same holds for $Y$). The tangent \b 
\ $TM$ restricted to $Y$ splits equivariantly as the direct sum of $TY 
$ and the normal bundle ${\mathcal N}(Y)$ which inherits a complex 
structure from the $S^1$-action. Let $x_i+m_{Y,i}\cdot z$, $d(Y)<i\leq 
m$, denote the equivariant roots of ${\mathcal N}(Y)$ (to lighten the 
notation we suppress the dependence of $x_i$ on $Y$). 

On $Y$ we choose the orientation which is compatible with the orientation of 
$M$ and the complex normal bundle ${\mathcal N}(Y)$. Let $\pm x_1,\ldots ,\pm x_{d(Y)}$
 denote a set of roots of $TY$ such that $x_1\cdot 
\ldots 
\cdot x_{d(Y)}$ is equal to the Euler class of the oriented \vb \ 
$TY$ and let $m_{Y,i}=0$ for $i\leq d(Y)$. Note that $\pm 
(x_i+m_{Y,i}\cdot z)$ are the equivariant roots of $TM_{\vert Y}$ as introduced in Section \ref{equiroots}.

Recall that the $Spin^c$-structure induces a complex line bundle over 
$M$. Let $l_Y$ be its weight at $Y$. The equivariant roots of $V$ at 
$Y$ shall be denoted by $v_1+s_{Y,1}\cdot z,\ldots ,v_s+s_{Y,s}\cdot 
z$. Next let 
$$I_Y:=\frac 1 2(\sum
_js_{Y,j}^2-\sum _im_{Y,i}^2),$$
which is an integer since we are looking at the two-fold action. 
Finally, define $n(V_{\vert Y}):=\dim _\C (V_0)$, where $V_0:=(V_{\vert Y})^{S^1}$ denotes 
the subbundle of $V_{\vert Y}$ which is fixed under the $S^1$-action. 

We are now in the position to state the main
 result of this section which we use in the following section to derive a lower bound
  for the first Pontrjagin class of a cohomology $\C P^m$ with $Pin(2)$-action.

\bigskip
\begin{proposition}\label{local data} For $q=e^{2\pi i\cdot 
\tau }$, $\tau \in \oha $, and $\lambda _0=e^{2\pi i\cdot z_0}$ a
 topological generator of $S^1\subset \C $ the series $\nu 
_Y(\tau ,z_0):=\tilde \nu _Y(q,\lambda _0)$ converges to a meromorphic function on $\oha 
\times \C $ also denoted by $\nu _Y(\tau ,z)$. The sum $\sum _Y\nu _Y(\tau ,z)$ has no poles on
$z\in \R $. 
\begin{enumerate}
\item If $n(V_{\vert Y})>d(Y)$ then $\nu _Y(\tau ,z)$ and $\tilde \nu 
_Y(q,\lambda )$ vanish identically.
\item If $n(V_{\vert Y})=d(Y)$ then $\nu 
_Y(\tau ,z)$ is the product of a holomorphic function $e(z)$ and a Jacobi function 
$F_Y$ for $SL_2(\Z )\ltimes 
\Z ^2$ of index $I_Y$. For any fixed $\tau \in \oha $ 
the set of poles of $F_Y$ is contained in $\Q 
\cdot
\tau +\Q
$.\end{enumerate}
\end{proposition}
 
\noindent{\bf Proof:} We describe the local datum $\tilde \nu _Y(q,\lambda _0)$ in terms
 of the Weierstra\ss ' $\Phi $-function
and the equivariant roots. Recall that $\Phi (\tau 
,z)$ is a holomorphic function on $\oha \times \C $ defined by the 
normally convergent infinite product $$\phi (q, \lambda ):=(\lambda 
^{\frac 12}-\lambda ^{-\frac 12})\cdot \prod _{n\geq 1}\frac 
{(1-q^n\cdot \lambda )\cdot (1-q^{n}\cdot \lambda ^{-1})}{(1-q^n)^2}\in \C \lbrack 
\lambda^{\frac 1 2}, \lambda^{-\frac 1 2}\rbrack \lbrack \lbrack q\rbrack \rbrack ,$$ 
where $q=e^{2\pi i\cdot \tau }$ and $\lambda =e^{2\pi i\cdot z }$. For 
a topological generator $\lambda _0$ of $S^1$ the power series $\tilde 
\nu 
_Y(q,\lambda _0) \in \C \lbrack 
\lbrack q\rbrack \rbrack $ is given by
\begin{equation}\label{equ1}\tilde \nu _Y(q,\lambda _0)=\left \langle A(q,\lambda _0),
 \mu _Y\right \rangle \end{equation}
where $A(q,\lambda )\in H^*(Y;\C (\lambda)\lbrack \lbrack q\rbrack \rbrack)$ is defined 
as\footnote{Since we have passed to the two-fold action this is an expression in $\lambda $
 rather than $\lambda^{\frac 1 2}$.}
$$e^{\frac 1 2 \cdot (c-\sum _jv_j)}\cdot \lambda 
^{\frac 1 2 \cdot (l_Y-\sum _js_{Y,j})}\times $$
$$\times \prod _{m_{Y,i}=0} \frac 
{x_i} {\phi (q,e^{x_i})}\cdot \prod _{m_{Y,i}\neq 0}\frac 1 {\phi (q, 
e^{x_i}\cdot \lambda ^{m_{Y,i}})}\cdot \prod_{j=1}^s \phi (q, 
e^{v_j}\cdot \lambda ^{s_{Y,j}}).$$
Here $\mu 
_Y$ is the fundamental cycle of $Y$, the characteristic classes are expressed in 
terms of their formal roots and $\langle \quad , \quad \rangle 
$ denotes the pairing between cohomology and homology. Recall that $l_Y$ is the 
weight at $Y$ of the complex line bundle associated to the equivariant $Spin^c$-structure.
 To prove formula (\ref{equ1}) one computes the Chern character of ${\mathcal U}_V$ and applies 
the \LFF \ to $ind_{S^1}(\partial _{\c}\otimes {\mathcal U}_V)(\lambda 
_0)$ (for details cf. \cite{DTOP}).

Recall that $\phi (q,\lambda )$ converges normally to $\Phi (\tau 
,z)$. This implies that for fixed $\lambda _0=e^{2\pi i\cdot z_0}$ and 
any $q=e^{2\pi i\cdot 
\tau }$, $\tau \in \oha $, the series $A(q,\lambda _0)\in H^*(Y;\C )\lbrack \lbrack 
q\rbrack \rbrack $ converges to a well defined element ${\mathcal A}_{z_0}(\tau )$ in the 
cohomology of $Y$ with values in the ring of holomorphic functions on 
$\oha $. Moreover, there exists an element ${\mathcal A}(\tau ,z)$ in the cohomology 
of $Y$ with values in the ring ${\mathcal M}(\oha \times \C )$ of 
meromorphic functions on $\oha 
\times 
\C $ such that ${\mathcal A}(\tau ,z_0)={\mathcal A}_{z_0}(\tau )$ for any irrational
 $z_0\in \R $ (for details cf. \cite{DTOP}). Changing variables we conclude from formula (\ref{equ1})
 that the series $\nu _Y(\tau ,z )$ converges for any 
irrational real number $z$ to the meromorphic function
 $\langle {\mathcal A}(\tau ,z),\mu _Y\rangle $ (in the following also denoted by $\nu _Y(\tau ,z)$). 

Note that each coefficient of the $q$-power series $ind_{S^1}(\partial 
_{\c}\otimes {\mathcal U}_V)$, being a finite Laurent polynomial in 
$\lambda $, is holomorphic on $S^1\subset \C $. This implies that the 
sum $\sum _Y\nu 
_Y(\tau ,z)$ has no poles on $z\in \R $ (again cf. \cite{DTOP} for details). Next we consider the 
statements involving the dimension $n(V_{\vert Y})$ of $V_0$.

Ad (1): Assume $n(V_{\vert Y})>d(Y)$. Recall that the Weierstra\ss ' 
$\Phi $-function $\Phi (\tau ,z)$ has a simple zero in $z=0$ for any 
$\tau\in \oha $. This implies that ${\mathcal A}(\tau ,z)$ contains 
the Euler class $e(V_0)=\prod _{s_{Y,j}=0}v_j$ of $V_0$ as a factor. 
Since $n(V_{\vert Y})>d(Y)$ the function $\nu 
_Y(\tau ,z)=\langle {\mathcal A}(\tau ,z),\mu _Y\rangle $ vanishes for any irrational $z\in \R $. 
Being meromorphic this forces $\nu 
_Y(\tau ,z)$ (and also $\tilde 
\nu 
_Y(q,\lambda )$) to vanish identically.

Ad (2): Assume $n(V_{\vert Y})=d(Y)$. Then $e(V_0)$ is in the top 
degree of the cohomology of $Y$. The local datum $\nu 
_Y(\tau ,z 
_0)$ is equal to the product of the Euler number of $V_0$ and ${\mathcal A}_0(\tau, z_0)$, where
 ${\mathcal A}_0(\tau, z)=e^{\pi \cdot i(l_Y-\sum 
_js_{Y,j})\cdot z}\cdot F_Y(\tau ,z)$ and $$F_Y(\tau ,z)=\prod 
_{m_{Y,i}\neq 0}\frac 1 {\Phi (\tau ,m_{Y,i}\cdot z)}\cdot \prod _{s_{Y,j}\neq 0} \Phi 
(\tau ,s_{Y,j}\cdot z).$$
To see this consider $A(q,\lambda )$, recall that $\phi (q,\lambda )$ converges to
 $\Phi (\tau ,z)$ and note that $\phi (q ,e^x)$ has the form $x + O(x^3)$. Whereas 
 ${\mathcal A}(\tau ,z)$ depends on the 
equivariant roots ${\mathcal A}_0(\tau, z)$ only depends on the weights 
of $TM$ and $V$ at $Y$.

We proceed to identify $F_Y(\tau ,z)$ with a 
Jacobi function. The Weierstra\ss ' $\Phi $-function is a holomorphic Jacobi 
function for $SL_2(\Z )\ltimes \Z ^2$ of weight $-1$ and `index $\frac 
1 2$ with character' (cf. \cite{EiZa}). More precisely, $\Phi $ is holomorphic and
satisfies 
$$\Phi (\tau ,z+\alpha \cdot \tau +\beta )=\Phi (\tau ,z)\cdot e^{-\pi i\cdot (\alpha ^2\tau +2\alpha
\cdot z)}\cdot (-1)^{\alpha +\beta }$$
for $(\alpha ,\beta )\in \Z ^2$ and
$$\Phi (\frac {a\tau +b}{c\tau +d},\frac z {c\tau +d})=\Phi (\tau ,z)\cdot (c\tau +d)^{-1}\cdot
e^{\pi i \cdot \frac{c\cdot z^2}{c\tau +d}}$$
for $\left ( \begin{smallmatrix}a&b\\c&d\end{smallmatrix}\right )\in SL_2(\Z )$.
 In particular, if $s$ is even, then $\Phi 
(\tau ,s\cdot z)$ is a holomorphic Jacobi 
function of weight $-1$ and index $\frac {s^2}2$ (as defined in the beginning of
 this section). For fixed $\tau $ the divisor of $\Phi _{\tau }(z):=\Phi 
(\tau ,z)$ is equal to $\Z 
\tau +\Z
$, i.e. $\Phi _{\tau }$ has a simple zero in each lattice point. From these 
properties it follows that $F_Y(\tau, z)$ is a Jacobi function for $SL_2(\Z )\ltimes 
\Z ^2$ of index $I_Y$ with poles in $\Q \cdot \tau 
+\Q $. This completes the proof of the proposition. 
\proofend

\section{A lower bound for the first Pontrjagin class}\label{linear bound}
In this section we give a lower bound for the first Pontrjagin 
class of a cohomology complex projective space which supports a nice 
$Pin(2)$-action with fixed point. Let $M$ be a $2m$-dimensional 
manifold with $H^*(M;\Z )\cong H^*(\C P^m;\Z )$. We assume that 
$Pin(2)$ acts almost effectively on $M$ and trivially on cohomology. Also 
we assume that the action has a fixed point. As mentioned before the 
latter assumption is automatically satisfied if the Euler 
characteristic of $M$ is odd, i.e. if $m$ is even. 

Let $\gamma $ denote the complex line bundle with $c_1(\gamma )=x$, 
where $x$ is a fixed generator of $H^2(M;\Z )$. Since the $Pin(2)$-action
 is trivial on integral cohomology, we can lift the action uniquely to $\gamma $
  (this follows form \cite{HaYo}, cf. \cite{DAN}, Prop. 3.6). Let $a_Y$ 
denote the weight of $\gamma $ at a connected component $Y$ of 
$M^{S^1}$ for the induced $S^1$-action.

Next we recall some well known facts about $M$ and $\gamma $ which may 
be induced from the localization theorem in $K$-theory applied to the 
$S^1$-action and induced $\Z _p$-actions or by cohomological means 
(cf. for example \cite{Pe1}, Th. 2.8, \cite{Br} Ch. VII, \cite{Hs}): 
\begin{enumerate}
\item[(i)] The fixed point
manifold $M^{S^1}$ is a disjoint union $Y_0\cup
\ldots \cup Y_k$, where the integral cohomology ring of $Y_i$ is isomorphic to $\C P^{m_i}$ for
some $m_i$.
\item[(ii)] $\sum _{i=0}^k
(m_i+1)=m+1$.
\item[(iii)] The weights of $\gamma $ are distinct.
\end{enumerate}
Assume $a_{Y_0}=0$, i.e. assume that the $Pin(2)$-fixed point 
$pt$ is in $Y_0$. Let $V$ be the sum of $S^1$-equivariant complex line bundles over $M$ given by 
$$V:=d(Y_0)\cdot\gamma  +\sum _{i=1}^k(d(Y_i)+1)\cdot \gamma \otimes 
\lambda ^{-a_{Y_i}},$$
where $\lambda $ denotes the standard complex one-dimensional 
representation of $S^1$. We apply Proposition \ref{local data} to 
$ind_{S^1}(\partial _{\c }\otimes {\mathcal U}_V)$ to derive a lower 
bound for the first Pontrjagin class. 

\bigskip
\begin{proposition}\label{bound on p1} Let $M$ be a cohomology $\C P^m$ as above. 
If $p_1(M)=-n\cdot x^2$ then $n<m$.\end{proposition}

\noindent{\bf Proof:} We may assume that $n$ is non-negative. Using the 
Atiyah-Singer index theorem one computes that the non-equivariant 
index $ind(\partial _{\c }\otimes {\mathcal U}_V)$ does not vanish. We 
proceed to describe the equivariant index in terms of local data. For 
technical reasons we replace the $S^1$-action by its two-fold action. 
Recall that $d(Y_i)$ denotes half of the dimension of $Y_i$, i.e. 
$d(Y_i)=m_i$, and $n(V_{\vert Y_i})$ denotes the complex dimension of 
the subbundle $V_0$ of $V_{\vert Y_i}$ which is fixed under the 
$S^1$-action. Since the weights $a_{Y_0},\ldots ,a_{Y_k}$ of $\gamma $ 
are distinct and $a_{Y_0}=0$ we have $n(V_{\vert Y_i})=d(Y_i)+1$ if 
$i>0$ and $n(V_{\vert Y_0})=d(Y_0)$. 

Next consider the local datum $\nu 
_{Y_i}(\tau ,z)$ in the \LFF \ for $ind_{S^1}(\partial _{\c }\otimes 
{\mathcal U}_V)$. By Proposition \ref{local data} $\nu _{Y_i}(\tau ,z)$ 
vanishes for $i>0$, $\nu _{Y_0}(\tau ,z)$ is the product of a 
holomorphic function $e(z)$ and a meromorphic function $F_{Y_0}(\tau ,z)$ 
and $ind_{S^1}(\partial _{\c }\otimes {\mathcal U}_V)(\lambda )$ converges 
to 
$$\sum _i\nu
_{Y_i}(\tau ,z)=\nu 
_{Y_0}(\tau ,z)=e(z)\cdot F_{Y_0}(\tau ,z)$$
for $q=e^{2\pi i\cdot \tau }$ and any topological generator $\lambda 
=e^{2\pi i\cdot z }$. Moreover $F_{Y_0}(\tau ,z)$ is a Jacobi function 
for $SL_2(\Z )\ltimes \Z ^2$ of index 
$$I_{Y_0}:=\frac 1 2(\sum _{i=1}^k(d(Y_i)+1)\cdot a_{Y_i}^2 -\sum _{j=1}^m m_{Y_0,j}^2)$$
with poles in $\Q \tau +\Q $. Note that 
$F_{Y_0}(\tau ,z)$ and $e(z)$ cannot vanish identically, since the 
non-equivariant index $ind(\partial _{\c }\otimes {\mathcal U}_V)$ does not vanish.

It follows from the proof of Proposition 
\ref{local data} that $e(z)$ has no zeros on $z\in \R $. Since $\sum 
_i\nu _{Y_i}(\tau ,z)$ has no poles on $z\in \R $ the same holds for 
$F_{Y_0}(\tau ,z)$. 

Next consider the action of $A=\left 
(\begin{smallmatrix}a&b\\c&d\end{smallmatrix}\right )\in SL_2(\Z )$ on 
$\oha \times \C $ given by $A(\tau , z)=(\frac {a\cdot \tau +b}{c\cdot 
\tau +d},\frac z {c\tau +d})$. Note that the $SL_2(\Z )$-orbit of any 
element of $\{\tau \} \times (\Q 
\tau +\Q
)$ intersects with $\{\tau _0\}\times \Q$ for some $\tau _0\in \oha $. 
Since $F_{Y_0}(\tau ,z)$ has poles in $\Q \tau +\Q $ but no poles on 
$\R $ we conclude that $F_{Y_0}$ has no poles at all. Hence, $F_{Y_0}$ 
is a holomorphic Jacobi function of index $I_{Y_0}$ which does not 
vanish identically. Since a holomorphic Jacobi function of negative 
index must vanish identically (cf. \cite{EiZa}) the index $I_{Y_0}$ is 
non-negative, i.e. 
$$\sum_{j=1}^m m_{Y_0,j}^2\leq \sum _{i=1}^k(d(Y_i)+1)\cdot a_{Y_i}^2.$$
Let $Z$ be a connected fixed point component of $M^{S^1}$ such that $a_Z^2=\max _i \{a_{Y_i}^2\}$. By
equation $(\ast )$
$$\sum _{j=1}^m m_{Y_0,j}^2=\sum _{j=1}^m m_{Z,j}^2+n\cdot a_Z^2.$$
Hence,
$$\sum _{j=1}^m m_{Z,j}^2+n\cdot a_Z^2\leq \sum _{i=1}^k(d(Y_i)+1)\cdot a_{Y_i}^2.$$
This implies $n<m$ since $a_Z^2\geq a_{Y_i}^2$ and $\sum 
_{i>0}(d(Y_i)+1)\leq m$.\proofend
 
\section{Rigidity of Pontrjagin classes}\label{proof}
In this section we prove Theorem \ref{main}. In fact we show the following slightly more 
general result. 

\bigskip
\begin{theorem}\label{maingeneral} Let $M$ be a smooth cohomology $\C P^{m}$, $m<12$, which supports
 an almost effective smooth $Pin(2)$-action which is trivial on cohomology. If $m$ is odd 
assume in addition that the action has a fixed point. Then the total 
Pontrjagin class of $M$ is standard, i.e. $p(M)=(1+x^2)^{m+1}$. 
\end{theorem}

\bigskip
\begin{remarks}
\begin{enumerate}
\item The theorem above is slightly more general than Theorem \ref{main} since a cohomology $\C
P^m$ may have non-trivial fundamental group.
\item For $m$ odd there are $S^3$-actions on $\C P^m$ with fixed point free $Pin(2)$-action. The
homogeneous action on $S^3/S^1\cong \C P^1$ is an example.
\end{enumerate}
\end{remarks}

\noindent
{\bf Proof of Theorem \ref{maingeneral}:} First note that the 
$Pin(2)$-action on $M$ always has a fixed point (for $m$ even this is 
true since the Euler characteristic of $M$ is odd). We order $H^4(M;\Z 
)=\Z \cdot x^2$ by identifying $H^4(M;\Z )$ with the integers using 
$x^2\mapsto 1$. By Proposition \ref{bound on p1} the first Pontrjagin 
class satisfies $p_1(M)>-m\cdot x^2$. In \cite{DAN}, Th. 4.2, it was 
shown that $p_1(M)\leq (m+1)\cdot x^2$ using methods similar to those 
of Section \ref{jacobi}. For the convenience of the reader we sketch 
the argument below (see Th. \ref{upperbound}). Hence, 
\begin{equation}\label{equ2}-m\cdot x^2<p_1(M)\leq (m+1)\cdot x^2.\end{equation}
Next note that $p_1(M)$ is a cohomology invariant modulo $24$, i.e. 
$p_1(M)\equiv (m+1)\cdot x^2\bmod 24$. To see this choose a $\spc 
$-structure on $M$ with first Chern class $c=(m+1)\cdot x$ and let 
$V:=(\gamma 
-1)^{m-2}$. By the Atiyah-Singer index theorem the index of the 
(non-equivariant) $Spin^c$-Dirac operator twisted with $V$ is equal to 
$$\langle e^{\frac c 2}\cdot \hat {\mathcal A}(M)\cdot (e^x-1)^{m-2},\mu _M\rangle .$$
Note that $\hat {\mathcal A}(M)=1 -\frac {p_1(M)}{24}+$ terms of 
higher order. Let $b$ be the integer defined by $p_1(M)=b\cdot x^2$. 
Then the index takes the form $\frac b {24}-Q$, where $Q$ is a 
rational number which only depends on the cohomology ring of $M$. 
Since the index is an integer it follows that $\frac b {24}\equiv Q$ 
modulo the integers. The same computation shows for the standard 
complex projective space that $\frac {m+1} {24}\equiv Q\bmod 
\Z$. Hence, $b\equiv m+1 
\bmod 24$, i.e.
$p_1(M)\equiv (m+1)\cdot x^2$ modulo $24$. Since $m<12$ it follows 
from equation (\ref{equ2}) that the first Pontrjagin class is 
standard, i.e. 
$$p_1(M)=(m+1)\cdot x^2.$$
 As mentioned in the 
introduction it was shown in \cite{DAN} that the total Pontrjagin 
class of $M$ is standard if $p_1(M)$ is standard (see Th. 
\ref{upperbound} below). This completes the proof.\proofend 

\bigskip
\begin{theorem}
\label{upperbound} Let $M$ be a cohomology $\C P^m$ with nice $Pin(2)$-action. If $m$ is odd 
assume in addition that the action has a fixed point. Then the first 
Pontrjagin class satisfies $p_1(M)\leq (m+1)\cdot x^2$. Moreover 
$p(M)$ is standard if $p_1(M)=(m+1)\cdot x^2$. 
\end{theorem}

\noindent{\bf Proof:} For a detailed proof we refer to \cite{DAN}, Th. 4.2. 
Here is a sketch of the argument based on the following general 
vanishing result. Let $V$ be a $Pin(2)$-equivariant complex 
\vb 
\ and let $W$ be a $Pin(2)$-equivariant $2t$-dimensional 
$Spin$-\vb 
\ over $M$. Let $\pm 
(w_1+t_{Y,1}\cdot z),\ldots , \pm (w_t+t_{Y,t}\cdot z)$ denote the 
equivariant roots of $W$ restricted to a connected component $Y$ of 
$M^{S^1}$. The equivariant roots of $V$ and $TM$ shall be denoted as 
in Section \ref{jacobi}. Assume $p_1(V+W)=p_1(M)$. Using a 
spectral-sequence argument one shows that 
\begin{equation}\label{equ3}I:=\frac 1 2 (\sum _{j=1}^s s_{Y,j}^2 +\sum _{k=1}^t
t_{Y,k}^2-\sum _{i=1}^m m_{Y,i}^2)\end{equation} is independent of 
$Y$. Next consider the $q$-power series ${\mathcal U}_{V,W}\in 
K_{S^1}(M)\powser q$ of virtual $S^1$-equivariant \vb s defined by 
$${\mathcal U}_{V,W}:={\mathcal U}_V\otimes \triangle (\widetilde {W})\otimes \twoindex {n=1} \infty
\bigotimes \Lambda _{q^n}(\widetilde {W}\otimes _\R \C ).$$ Here $\triangle (W)$ denotes the full
complex spinor \b \ associated to the $Spin$-\vb \ $W$. We fix a $\spc 
$-structure on $M$ and lift the induced $S^1$-action. Next one shows 
(by arguments similar to the ones of the previous sections) that the 
equivariant index $ind_{S^1}(\partial _{\c }\otimes {\mathcal 
U}_{V,W})$ is equal to the product of a holomorphic function and a 
Jacobi function (for $\Gamma _0(2)\subset SL_2(\Z )$) of index $I$ if 
$p_1(V+W)=p_1(M)$ and $c_1(V)$ is equal to the first Chern class of 
the $\spc $-manifold $M$. As in the proof of the rigidity of elliptic 
genera one can show that the Jacobi function is in fact holomorphic. 
This implies that $ind_{S^1}(\partial 
_{\c }\otimes {\mathcal U}_{V,W})$ vanishes identically if $I$ is 
negative. 

Now assume $p_1(M)>(m+1)\cdot x^2$. We want to show a contradiction. 
To this end fix a $\spc $-structure on $M$ with first Chern class 
equal to $(m+1)\cdot x$. Let $V:=(m-1)\cdot \gamma +\gamma ^2$ and let 
$W:=(b-m-3)\cdot \gamma $, where $b>m+1$ is defined by $p_1(M)=b\cdot 
x^2$ (note that $b\geq m+3$ since $p_1(M)\equiv (m+1)\cdot x^2\bmod 
2$). We lift the $Pin(2)$-action to each line bundle occurring in $V$ 
and $W$. Note that $p_1(V+W)=p_1(M)$ and $c_1(V)$ is equal to the 
first Chern class of $M$. Since the weights of $V$ and $W$ at the 
$Pin(2)$-fixed point vanish it follows from equation (\ref{equ3}) that 
$I$ is negative. Hence, by the result above $ind_{S^1}(\partial _{\c 
}\otimes {\mathcal U}_{V,W})$ vanishes identically. In particular, the 
non-equivariant index vanishes. However, one computes with the help of 
the Atiyah-Singer index theorem that the series $ind(\partial _{\c 
}\otimes {\mathcal U}_{V,W})$ does not vanish. This contradicts the 
assumption on $p_1(M)$. Thus $p_1(M)\leq (m+1)\cdot x^2$. 

Next assume $p_1(M)=(m+1)\cdot x^2$. We want to show that $p(M)$ is 
standard. To this end let $V_k:=\gamma 
^2+(m-3-2k)\cdot 
\gamma $, $W_k:=(2k)\cdot 
\gamma $, $k\in \{0,\ldots ,\lbrack \frac {m-3} 2\rbrack \}$, and choose a $Spin^c$-structure on $M$ with 
first Chern class equal to $c_1(V_k)$. Note that 
$p_1(V_k+W_k-TM)=0$. Again we lift the $Pin(2)$-action to each line 
bundle occurring in $V_k$ and $W_k$ and conclude from equation 
(\ref{equ3}) that $I$ is negative. By the result above 
$ind_{S^1}(\partial 
_{\c }\otimes {\mathcal U}_{V,W})$ vanishes identically. In 
particular, the constant term in the $q$-power series is zero, i.e. 
$$\left \langle \hat {\mathcal A}(M)\cdot (e^x-e^{-x})\cdot (e^{\frac x 2} -e^{-\frac x
2})^{m-3-2k}\cdot (e^{\frac x 2}+e^{-\frac x 2})^{2k}, \mu _M
\right \rangle =0$$
for $k\in \{0,\ldots ,\lbrack \frac {m-3} 2\rbrack \}$. These 
relations together with the signature theorem completely determine 
$\hat {\mathcal A}(M)$ and therefore determine the total Pontrjagin 
class $p(M)$. Since all these relations also hold true for $\C P^m$ we 
conclude that $p(M)=(1+x^2)^{m+1}$. Hence $p(M)$ is standard if 
$p_1(M)$ is standard.\proofend 

\setlength{\parindent}{0pt}

\vskip3cm
\noindent
Anand Dessai\\e-mail: dessai@math.uni-augsburg.de\\ 
http://www.math.uni-augsburg.de/geo/dessai/homepage.html\\ Department 
of Mathematics, University of Augsburg, D-86159 Augsburg 

\end{document}